\documentclass[letterpaper, 12 pt]{article}
\usepackage{amsthm, amsmath, amssymb, amsfonts, doublespace}

\newcommand{\setleftmargin}[1]{
        \addtolength{\textwidth}{\oddsidemargin}
        \addtolength{\textwidth}{1in}
        \addtolength{\textwidth}{-#1}
        \setlength{\oddsidemargin}{-1in}
        \addtolength{\oddsidemargin}{#1}
        \setlength{\evensidemargin}{\oddsidemargin}
}
\newcommand{\setrightmargin}[1]{
        \setlength{\textwidth}{8.5in}
        \addtolength{\textwidth}{-\oddsidemargin}
        \addtolength{\textwidth}{-1in}
        \addtolength{\textwidth}{-#1}
}
\setleftmargin{1 in}
\setrightmargin{1 in}

\newcommand{\CC}{\mathbb{C}}
\newcommand{\PP}{\mathbb{P}}

\newcommand{\Cn}{\CC^n}
\newcommand{\Cnd}{(\Cn)^{\vee}}

\newcommand{\into}{\hookrightarrow}
\newcommand{\onto}{\twoheadrightarrow}

\DeclareMathOperator{\In}{In}

\DeclareMathOperator{\nbc}{bc}

\DeclareMathOperator{\Spec}{Spec}

\DeclareMathOperator{\SR}{\mathrm{SR}}

\newcommand{\A}{\mathcal{A}}

\newtheorem{conj}{Conjecture}
\newtheorem{Theorem}[conj]{Theorem}
\newtheorem{prop}[conj]{Proposition}
\newtheorem{lemma}[conj]{Lemma}

\newtheorem{remark}[conj]{Remark}

\begin{document}
\spacing{1.1}

\noindent
{\Large \bf A Broken Circuit Ring}
\bigskip\\
{\bf Nicholas Proudfoot}\footnote{Partially supported
by the Clay Mathematics Institute Liftoff Program}\\
Department of Mathematics, University of Texas,
Austin, TX 78712\smallskip \\
{\bf David Speyer} \\
Department of Mathematics, University of California,
Berkeley, CA 94720
\bigskip
{\small
\begin{quote}
\noindent {\em Abstract.}
Given a matroid $M$ represented by a linear subspace $L\subset\Cn$
(equivalently by an arrangement of $n$ hyperplanes in $L$),
we define a graded ring $R(L)$ which degenerates to the Stanley-Reisner
ring of the broken circuit complex for any choice of ordering of the ground set.
In particular, $R(L)$ is Cohen-Macaulay, and may be used to compute
the $h$-vector of the broken circuit complex of $M$.
We give a geometric interpretation of $\Spec R(L)$, as well
as a stratification indexed by the flats of $M$.
\end{quote}
}

\section{Introduction}\label{intro}
Consider a vector space with basis $\Cn = \CC\{e_1,\ldots,e_n\}$, and
its dual $\Cnd = \CC\{x_1,\ldots,x_n\}$.
Let $L \subset \CC^n$ be a linear subspace of dimension $d$. 
We define a matroid $M(L)$ on the ground set 
$[n] := \{ 1, \ldots, n \}$ by declaring $I \subset [n]$ to be independent if and only 
if the composition $\CC\{x_i \mid i \in I\}\into \Cnd \onto L^\vee$ is injective. 
Recall that a minimal dependent subset $C \subset [n]$ is called a \emph{circuit}; 
in this case there exist scalars $\{a_c\mid c\in C\}$, unique up to scaling, 
such that $\sum_{C} a_c x_c$ vanishes on $L$.  Conversely, the support
of every linear form that vanishes on $L$ contains a circuit.

The central object of study in this paper will be the ring $R(L)$ generated
by the inverses of the restrictions of the linear functionals $\{x_1,\ldots,x_n\}$
to $L$.  More formally, let
$$\CC[x,y] := \CC[x_1,y_1,\ldots,x_n,y_n]/\langle x_iy_i - 1\rangle,$$
and let $\CC[x]$ and $\CC[y]$ denote the polynomial subrings generated
by the $x$ and $y$ variables, respectively.
Let $\CC[L]$ denote the ring of functions on $L$,
which is a quotient of $\CC[x]$ by the ideal generated by the linear forms
$\big\{\sum_C a_c x_c\mid C\text{ a circuit}\big\}$.
We now set $$R(L) := \Big(\CC[L] \otimes_{\CC[x]} \CC[x,y]\Big) \cap\CC[y].$$
Geometrically, $\Spec R(L)$ is a subscheme of $\Spec\CC[y]$, which we will
identify with $\Cnd$.  Using the isomorphism between $\Cn$ and $\Cnd$ provided
by the dual bases, $\Spec R(L)$ may be
obtained by intersecting $L$ with the torus $(\CC^*)^n$, 
applying the involution $t \mapsto t^{-1}$ on the torus, and taking the closure inside 
of $\CC^n$. If $C$ is any circuit of 
$M(L)$ with $\sum_{c \in C} a_c x_c$ vanishing on $L$, then we have 
the relation $$f_C := \sum_{c \in C} a_c \!\!\prod_{c'\in C\setminus\{c\}}\!\!y_{c'}=0
\,\,\,\text{  in  }\,\,\, R(L).$$  
Our main result (Theorem \ref{degen}) will be that the elements $\{f_C\mid C\text{ a circuit}\}$ 
are a universal Gr\"obner basis for $R(L)$, hence this ring degenerates to the Stanley-Reisner ring
of the broken circuit complex of $M(L)$ for {\em any} choice of ordering
of the ground set $[n]$.
It follows that $R(L)$ is a Cohen-Macaulay ring of dimension $d$, and that the quotient
of $R(\A)$ by a minimal linear system of parameters has Hilbert series 
equal to the $h$-polynomial
of the broken circuit complex.
In Proposition \ref{free} we identify a natural choice of linear parameters for $R(L)$.

The Hilbert series of $R(L)$ has already been computed by Terao \cite{Te},
using different methods.
The main novelty of our paper lies in our geometric approach, 
and our interpretation of $R(L)$ as a deformation of another well-known ring.
The ring $R(L)$ also appears as a cohomology ring in \cite{PW}, and as the homogeneous
coordinate ring of a projective variety in \cite[3.1]{Lo}.

\paragraph{\bf Acknowledgment.}  Both authors would like to thank Ed Swartz
for useful discussions.

\begin{section}{The broken circuit complex}\label{bcsec}
Choose an ordering $w$ of $[n]$. We define a \emph{broken circuit of $M(L)$ with respect to 
$w$} to be a set of the form $C \setminus \{ c \}$, where $C$ is a circuit of 
$M(L)$ and $c$ the $w$-minimal element of $C$. 
We define the \emph{broken circuit complex} $\nbc_w(L)$ to the simplicial complex on the ground
set $[n]$ whose faces are those subsets of $[n]$ that do not contain any broken circuit.
Note that all of the singletons will be faces of $\nbc_w(L)$ if and only if $M(L)$
has no parallel pairs, and the empty set will be a face if and only if $M(L)$ has no loops.
We will not need to assume that either of these conditions holds.

Consider the $f$-vector $(f_0, \ldots, f_d)$ of $\nbc_w(L)$,
where $f_i$ is the number of faces of order $i$.  Then 
$f_i$ is equal to the rank of $H^i(A(L))$, where
$A(L) = L \setminus \bigcup_{i=1}^n \{ x_i =0\}$ is the complement of the restriction
of the coordinate arrangement from $\Cn$ to $L$ (see for example \cite{OT}).  
In particular, the $f$-vector
of $\nbc_w(L)$ is independent of the ordering $w$.
The $h$-vector $(h_0,\ldots,h_{d-1})$ of $\nbc_w(L)$
is defined by the formula $\sum h_i z^i = \sum f_i z^i (1-z)^{d-i}$.

The \emph{Stanley-Reisner ring} $\SR(\Delta)$ of a simplicial complex 
$\Delta$ on the ground set $[n]$ is defined to be the quotient of 
$\CC[e_1, \ldots, e_n]$ by the ideal generated by the monomials
$\prod_{i \in N} e_i$, where $N$ ranges over the nonfaces of $\Delta$. 
The complex $\nbc_w(L)$ is shellable of dimension $d-1$ \cite{Bj}, which 
implies that $\Spec \SR(\nbc_w(L))$ is Cohen-Macaulay and pure of dimension $d$. 
Let $\CC[L^{\vee}]$ denote the ring of functions 
on $L^{\vee} = \Cnd/L^{\perp}$, which we may think of as the symmetric algebra on $L$. 
The inclusion of $L$ into $\Cn$ induces an inclusion of $\CC[L^{\vee}]$
into $\CC[e_1,\ldots,e_n]$, which
makes $\SR(\nbc_w(L))$ into an $\CC[L^{\vee}]$-algebra.
Let $\SR_0(\nbc_w(L)) = \SR(\nbc_w(L)) \otimes_{\CC[L^{\vee}]} \CC$,
where each linear function on $L^{\vee}$ acts on $\CC$ by $0$.
The following proposition asserts that $L$ constitutes a linear
system of parameters (l.s.o.p.) for $\SR(\nbc_w(L))$.

\begin{prop}\label{lsop}
The Stanley-Reisner ring $\SR(\nbc_w(L))$ is a free $\CC[L^{\vee}]$-module,
and the ring $\SR_0(\nbc_w(L))$ is zero-dimensional with Hilbert series $\sum h_i z^i$.
\end{prop}

\begin{proof}
By \cite[5.9]{St}, it is enough to prove that  
$\SR_0(\nbc_w(L))$ is a zero-dimensional ring. 
Let $\pi$ denote the composition  
$\Spec \SR(\nbc_w(L)) \into \Cnd \onto L^{\vee}$.
The variety $\Spec \SR(\nbc_w(L))$ is a union of coordinate subspaces, one for each 
face of $\nbc_w(L)$. Let $F$ be such a face, with vertices $(v_1, \ldots, v_{|F|})$. 
The broken circuit complex is a subcomplex of the matroid complex,
hence $(v_1, \ldots, v_{|F|})$ is an independent set, 
which implies that $\pi$ maps the corresponding coordinate subspace 
injectively to $L^{\vee}$. Thus $\pi^{-1}(0) = \Spec \SR_0(\nbc_w(L))$ is supported at the
origin, and we are done.
\end{proof}
\end{section}

\section{A degeneration of \boldmath$R(L)$}\label{mainsection}
In this section we show that $R(L)$ degenerates flatly
to the Stanley-Reisner ring $\SR(\nbc_w(L))$ for any choice of $w$.  

\begin{lemma}\label{dimdeg}
The spaces $\Spec R(L)$ and $\Spec \SR(\nbc_w(L))$ are both pure 
$d$-dimensional homogeneous varieties of degree $t_{M(L)}(1,0)$, where $t_M(w,z)$ 
is the Tutte polynomial of $M$.
\end{lemma}

\begin{proof}
The broken circuit complex is pure of dimension $d-1$, hence
$\Spec \SR(\nbc_w(L))$ is union of $d$-dimensional coordinate subspaces of $\Cnd$.
Its degree is the number of facets of $\nbc_w(L)$,
which is equal to $\sum h_i = t_{M(L)}(1,0)$ \cite{Bj}.

The variety $\Spec R(L)$ is equal to the closure inside of $\Cnd\cong\Cn$
of $L \cap (\CC^*)^n$, and is therefore $d$ dimensional. 
We will now show that $\deg \Spec R(L)$ obeys the same recurrence as
$t_{M(L)}(1,0)$. 
First, suppose that $i \in [n]$ is a loop of
$M(L)$. Then $L$ lies in a coordinate subspace of $\CC^n$, $L \cap
(\CC^*)^n$ is empty, and $\Spec R(L)$ is thus empty and has degree
$0$. In this case, we also have $t_{M(L)}(1,0)=0$. Next, suppose that
$i$ is a coloop of $M(L)$. Then $L$ is invariant under translation by
$e_i$, and $\Spec R(L)$ is similarly invariant under
translation by $x_i$. Write $L/i$ for the
quotient of $L$ by this translation, so that $\Spec R(L)=\Spec R(L/i)
\times \CC$ and $\deg \Spec R(L)=\deg \Spec R(L/i)$. It is clear that
$M(L/i)=M(L)/i$, and indeed $t_M(1,0)=t_{M/i}(1,0)$ when $i$ is a
coloop.

Now consider the case where $i$ is neither a loop nor a coloop, hence we
have $$t_{M(L)}(1,0)=t_{M(L)/i}(1,0) + t_{M(L) \setminus i}(1,0).$$ 
In this case, we may apply the following theorem.

\begin{Theorem}\label{kmy}{\em \cite[2.2]{KMY}}
Let $X$ be a homogeneous irreducible subvariety of $\CC^n = H \oplus \ell$,
with $H$ a hyperplane and $\ell$ a line such that $X$ is not invariant under
translation in the $\ell$ direction.
Let $X_1$ be the closure of the projection along $\ell$ of
$X$ to $H$, and let $X_2$ be the flat limit in $H \times \PP^1$ of
$X \cap (H \times \{ t \})$ as $t \to \infty$. Then $X$ has a flat
degeneration to a scheme supported on $(X_1 \times \{ 0 \}) \cup (X_2 \times \ell)$. In 
particular, $\deg X \geq \deg X_1+\deg X_2$, with equality if the projection $X \to X_1$ is 
generically one to one.
\end{Theorem}

Let $X=\Spec R(L)$, $\ell = \CC x_i$, and $H = \CC\{x_j\mid j\neq i\}$.
Then in the notation of Theorem \ref{kmy}, we have $X_1=\Spec R(L \setminus i)$, 
where $L \setminus i$ is the projection of $L$ onto $H$, 
and $X_2=\Spec R(L/i)$. The projection of $\Spec R(L)$ onto $H$ is one to one 
because the corresponding projection of $L$ in the $x_i$ 
direction is one to one. Thus the degree of $\Spec R(L)$ is additive.
\end{proof}

We are now ready to prove our main theorem, which asserts that $R(L)$
degenerates flatly to $\SR(\nbc_w(L))$ for any choice of $w$.

\begin{Theorem}\label{degen}
The set $\big\{f_C\mid C\text{ a circuit of $M(L)$}\big\}$ is a universal Gr\"obner basis
for $R(L)$.
Given any ordering $w$ of $[n]$, with the induced term order on $\CC[y]$,
we have $\In_w R(L)=\SR(\nbc_w(L))$.
\end{Theorem}

\begin{proof}
Suppose given an ordering $w$ of $[n]$ and a circuit $C$ of $M(L)$.
Let $c_0$ denote the $w$ minimal element of $C$, so that
$\prod_{c' \in C \setminus \{ c_0 \}} y_{c'}$ is the leading term of $f_C$ with respect to 
$w$.  Every monomial of this form vanishes in $\In_w R(L)$, hence we deduce that $\Spec 
\In_w(R(L))$ is a subscheme of $\Spec \SR(\nbc_w(L))$. However, 
Lemma \ref{dimdeg} tells us that these two schemes have 
the same dimension and degree, and $\Spec \SR(\nbc_w(L))$ is reduced. 
Thus they are equal.

Let $R$ be the quotient ring of $\CC[y]$ generated by the polynomials $\{f_C\}$. 
It is clear that $\In_w \Spec(R(L)) \subseteq \In_w 
\Spec R \subseteq \Spec \SR(\nbc_w(L))$. Since the two ends of this chain are equal, we have 
$\In_w R=\In_w R(L)$, and thus $R$ and $R(L)$ have the same Hilbert series. As $R(L)$ is a 
quotient ring of $R$, $R=R(L)$.
\end{proof}

\section{A stratification of \boldmath$\Spec R(L)$}
Let $I$ be a subset of $[n]$.  The {\em rank} of $I$ is defined to be the
cardinality of the largest independent subset of $I$. 
If any strict superset of $I$ has strictly greater rank, then $I$ is called a {\em flat}
of $M(L)$.  If $I$ is a flat, let 
$L_I\subset\CC^I$ be the projection of $L$ onto the coordinate subspace
$\CC^I\subset\Cn$, and let
$L^I\subset\CC^{I^c}$ be the intersection of $L$ with the complimentary
coordinate subspace $\CC^{I^c}$.
The matroid $M(L_I)$ is called the {\em localization of $M(L)$ at $I$},
while $M(L^I)$ is called the {\em deletion of $I$ from $M(L)$}.

For any $I \subset [n]$, let $U_I = \{y\in\Cnd\mid y_i = 0\iff i\notin I\}$,
and let $A_I = \Spec R(L) \cap U_I$.

\begin{prop}\label{strat}
The variety $A_I$ is nonempty if and only if $I$ is a flat of $M(L)$. 
If nonempty, $A_I$ is isomorphic to $A(L_I) = L_I\setminus\bigcup_{i\in I}\{y_i=0\}$.
\end{prop}

\begin{proof}
First suppose that $I$ is not a flat of $M(L)$. Then there exists some circuit $C$ of 
$M(L)$ and element $c_0\in C$ such that $C \cap I=C \setminus \{ c_0 \}$.
On one hand, the polynomial $f_C = \sum_{c \in C} 
a_c \prod_{c' \in C \setminus \{ c \}} y_{c'}$ vanishes on $A_I$.
On the other hand, $f_C$ has a unique nonzero term
$\prod_{c \in C \setminus \{ c_0 \}} y_{c'}$ on $U_I$, 
and therefore cannot vanish on this set.  Hence $A_I$ must be empty.

Now suppose that $I$ is a flat.
If $I=[n]$, then we are simply
repeating the observation that $\Spec R(L) \cap (\CC^*)^n \cong L \cap (\CC^*)^n = A(L)$.
In the general case, Theorem \ref{degen} tells us that $\Spec R(L)$ is cut out of 
$\Cnd$ by the polynomials $f_C$, so we need to understand the restrictions
of these polynomials to the set $U_I$.
If $C$ is not contained in $I$, then $C \setminus I$ has size at least $2$,
and therefore $f_C$ vanishes on $U_I$.  Thus we may restrict our attention
to those circuits that are contained in $I$.
Proposition \ref{strat} then
follows from the fact that the circuits of $M(L_I)$ are precisely
the circuits of $M(L)$ that are supported on $I$.
\end{proof}

\begin{remark}
The stratification of $\Spec R(L)$ given by Proposition \ref{strat}
is analogous to the standard stratification of $L$ into pieces isomorphic
to $A(L^I)$, again ranging over all flats of $M(L)$.
\end{remark}

The identification of $e_i$ with $y_i$ makes $R(L)$ into an algebra over $\CC[L^{\vee}]$.
We conclude by showing that, as in Proposition \ref{lsop},
$L$ provides a natural linear system of parameters for $R(L)$.

\begin{prop}\label{free}
The ring $R(L)$ is a free module over $\CC[L^{\vee}]$.
The zero dimensional quotient $R_0(L) := R(L)\otimes_{\CC[L^{\vee}]} \CC$
has Hilbert series $\sum h_i z^i$.
\end{prop}

\begin{proof}
The fact that $R(L)$ is Cohen-Macaulay follows from Theorem \ref{degen},
which asserts that it is a deformation of the Cohen-Macaulay ring 
$\SR(\nbc_w(L))$.  Furthermore, Theorem \ref{degen} tells us that
any quotient of $R(L)$ by $d$ generic parameters has the same Hilbert series
of $SR_0(bc_w(L))$.  Therefore, as in Proposition \ref{lsop},
we let $\pi$ denote the composition
$\Spec R(L)\into \Cnd \onto L^{\vee}$, 
and observe that it is enough to show that $\pi^{-1}(0)$ is supported at the origin.

Let $I \subset [n]$ and suppose that 
$y=(y_1, \ldots, y_n) \in A_I = \Spec R(L) \cap U_I$. 
By Proposition \ref{strat}, $A_I$ is obtained from $A(L_I)$
by applying the inversion involution of $(\CC^*)^I$, hence there
exists $x_I\in A(L_I)\subset L_I$ such that $x_i = y_i^{-1}$ for all $i\in I$.
Extend $x_I$ to an element $x\in L$.  Then $\langle x,y \rangle = \sum x_iy_i = |I|$,
hence if $y$ projects trivially onto $L^{\vee}$, we must have $I = \emptyset$.
\end{proof}

\begin{remark}
It is natural to ask the question of whether $R_0(L)$ has a $g$-element; that is 
an element $g\in R(L)$ in degree 1 such that the multiplication map
$g^{r-2i}:R_0(L)_i\to R_0(L)_{r-i}$ is injective for all $i< r/2$,
where $r$ is the top nonzero degree of $R_0(L)$.  This property is known
to fail for the ring $SR_0(bc_w(L))$ {\em\cite[\S 5]{Sw}}, but the inequalities that it would
imply for the $h$-numbers are not known to be either true or false.
In fact, the ring $R_0(L)$ fares no better than its degeneration;
Swartz's counterexample to the $g$-theorem for $SR_0(bc_w(L))$ is also a counterexample
for $R_0(L)$.
\end{remark}

\begin{remark}
All of the constructions and results in this paper generalize
to arbitrary fields with the exception of Proposition \ref{free}, which uses in an essential
manner the fact that $\CC$ has characteristic zero.
\end{remark}

\footnotesize{

}

\end{document}